%% file: main.tex
\documentclass[12pt]{amsart}
\usepackage[utf8]{inputenc}
\usepackage{graphicx}
\graphicspath{ {images/} }
\usepackage{enumerate}
\input{header}

\input{help}

\newcommand{\Qc}{\mathcal{Q}}
\title{Counting rational points close to $p$-adic integers and applications in Diophantine approximation}
\author{Benjamin Ward}
\address{Department of Mathematics, University of York, Heslington, York, YO10
5DD, United Kingdom}
\email{bw744@york.ac.uk}
\date{January, 2020}

\begin{document}
\begin{abstract}
We find upper and lower bounds on the number of rational points that are $\psi$-approximations of some $n$-dimensional $p$-adic integer. Lattice point counting techniques are used to find the upper bound result, and a Pigeon-hole principle style argument is used to find the lower bound result. We use these results to find the Hausdorff dimension for the set of $p$-adic weighted simultaneously approximable points intersected with $p$-adic coordinate hyperplanes. For the lower bound result we show that the set of rational points that $\tau$-approximate a $p$-adic integer form a set of resonant points that can be used to construct a local ubiquitous system of rectangles.
\end{abstract} 
\maketitle
\section{Introduction}
The study of rational points on algebraic varieties, usually called Diophantine geometry, has a wide variety of applications in many areas of mathematics. A variation of this is the study of rational points that lie close to such algebraic varieties. In the setting of $\R^{n}$ there has been many results of this type, including counts on the number of rational points close to curves \cite{BDV07,VV06,S16,S17,H96} and manifolds \cite{B12, BVVZ17, H20, HL18}. In the $p$-adic setting less is known. In \cite{BB20, BBS20} a bound on the number of rational points that lie on the curve $\Cf=\{(x,x^{2}, \dots, x^{n}): x \in \Zp \}$ were found, but as yet no other results are available. In this paper we provide an upper and lower bound on the number of rational points within a small neighbourhood of a $p$-adic integer. Such result allows us to find bounds on the number of rational points close to $p$-adic coordinate hyperplanes. \par

Fix a prime number $p \in \N$ and let $|.|_{p}$ denote the $p$-adic norm. Define the set of $p$-adic numbers $\Qp$ as the completion of $\Q$ with respect to the $p$-adic norm. Denote by $\Zp:=\{x \in \Qp: |x|_{p} \leq 1\}$ the ring of $p$-adic integers. Let $\bx \in \Zp^{n}$, $N \in \N$, and $\Psi=(\psi_{1}, \dots , \psi_{n})$ be an $n$-tuple of approximation functions of the form $\psi_{i}: \N \to \R_{+}$, with $\psi_{i}(q) \to 0$ as $q \to \infty$ for each $1 \leq i \leq n$. We provide bounds on the cardinality of the set
\begin{equation*}
\Qc(\bx,\Psi, N):=\left\{ (q_{0},q_{1}, \dots , q_{n}) \in  \Z^{n+1} : \begin{array}{c}
0<q_{0} \leq N, \\
\max_{1 \leq i \leq n}|q_{i}| \leq N, \\
\end{array}
 \, \left|q_{0} x_{i}-q_{i} \right|_{p} < \psi_{i}(N), \, 1 \leq i \leq n \right\}.
\end{equation*}
If the approximation functions $\psi_{i}$ are of the form $\psi_{i}(q)=q^{-\tau_{i}}$ for some vector $\bt=(\tau_{1}, \dots, \tau_{n}) \in \R^{n}_{>0}$ we will use the notation $\Qc(\bx, \bt, N)$. Note that to get a result for general $\bx \in \Zp^{n}$ we must apply some conditions. For example, if $\bx \in \Q^{n}$ then for sufficiently large $N \in \N$ we have that $\#\Qc(\bx,\Psi,N) \asymp N^{2}$ for any $\Psi$. Here $a \asymp b$ means there exists constants $c_{1}, c_{2} \in \R_{>0}$ such that $c_{1}b \leq a \leq c_{2}b$.  Conversely, if $\bx$ is badly approximable each approximation function satisfies $\psi_{i}(q)<q^{-1-\frac{1}{n}-\epsilon}$ for some $\epsilon>0$, then $\#\Qc(\bx,\Psi,N) \ll 1$. In order to obtain good bounds on the cardinality of $\Qc(\bx,\Psi,N)$ we use the Diophantine exponent $\tau(\bx)$ defined as 
\begin{equation*}
\tau(\bx):= \sup \left\{ \sum_{i=1}^{n}\tau_{i}: |q_{0}x_{i}-q_{i}|_{p} < Q^{-\tau_{i}}, \, \, \text{ for i.m. } \, Q \in \N \text{ with } |q_{i}| \leq Q \right\}.
\end{equation*}  
By a Theorem of Mahler \cite{S67} we have that for all $x \in \Zp$, $\tau(x) \geq 2$. Further, by a result of Jarnik \cite{J45} we have that $\tau(\bx) = n+1$ for almost all $\bx \in \Zp^{n}$, with respect to the $n$-dimensional Haar measure $\mu_{p,n}$ on $\Qp$, normalised by $\mu_{p,n}(\Zp^{n})=1$. \par 
We have the following result on the cardinality of $\Qc(x,\psi,N)$ for general $x \in \Zp$.
\begin{lemma} \label{lemma1}
Let $x \in \Zp$ with Diophantine exponent $\tau(x)$ and let $\psi(q)=q^{-\tau}$ for some $\tau \in \R_{+}$ with $\max \{1, \tau(x)-1\} < \tau < \tau(x)$. Then for any $\epsilon>0$ there exists sufficiently large $N_{0} \in \N$ such that for all $N \geq N_{0}$
\begin{equation*}
\#\Qc(x, \tau, N) \leq N^{\tau(x)-\tau+\epsilon}.
\end{equation*}
\end{lemma}
Note by our previous remark on the Diophantine exponent that for almost all $x \in \Zp$ we have $\tau(x)=2$, so the above lemma reads that for $\psi(q)=q^{-\tau}$ with $1<\tau<2$, then for almost all $x \in \Zp$
\begin{equation*}
\#\Qc(x,\psi, N) \leq N^{2-\tau+\epsilon}.
\end{equation*}
While Lemma \ref{lemma1} gives us an upper bound for all $x \in \Zp$, provided the approximation function $\psi$ is 'close' to the function related to the Diophantine exponent, the bound given has an extra $Q^{\epsilon}$ term, which we believe is unnecessary. The following theorem offers an improvement in this respect.
\begin{theorem} \label{theorem1}
Let $\bx \in \Zp^{n}$ and suppose that $\tau(\bx) = n+1$. Let $\Psi$ be an $n$-tuple of approximation functions with each
\begin{equation*}
 q^{-1-\frac{1}{n}+\epsilon} < \psi_{i}(q) < q^{-1}, \quad 1 \leq i \leq n,
\end{equation*} 
 for some $\epsilon>0$. Then there exists $N_{0} \in \N$ such that for all $N \geq N_{0}$,
\begin{equation*}
\#\Qc(\bx,\Psi,N) \leq C_{1}N^{n+1}\prod_{i=1}^{n}\psi_{i}(N),
\end{equation*}
where
\begin{equation*}
 C_{1}=\max\left\{ 3(6\sqrt{n})^{n}, \frac{(n+2)! \pi^{n/2}\sqrt{n}^{n+1}}{\Gamma\left(\frac{n}{2}+1 \right)}\right\}.
 \end{equation*}
\end{theorem}
Again, as with Lemma \ref{lemma1}, we can deduce that the above upper bound is true for almost all $\bx \in \Zp^{n}$. This type of result has already been proven in the real case (see Lemma 6.1 of \cite{BHV20}). In the case where the approximation functions are of the form $\psi_{i}(q)=q^{-\tau_{i}}$ then the theorem reads: if 
\begin{equation*}
\sum_{i=1}^{n} \tau_{i}<n+1, \quad \text{ and } \quad \tau_{i}>1,
\end{equation*}
then for any $\bx \in \Zp^{n}$ with $\tau(\bx)=n+1$,
\begin{equation*}
\#\Qc(\bx,\bt,N) \leq C_{1}N^{n+1-\sum_{i=1}^{n}\tau_{i}}.
\end{equation*}
\par 
 Lastly, we have the following lemma which provides a complimentary lower bound to the previous two results.
\begin{lemma} \label{lemma2}
 Let $\bx \in \Zp^{n}$ and 
\begin{equation*}
\sum_{i=1}^{n}\tau_{i}<n+1, \quad \text{and} \quad \tau_{i}>1
\end{equation*}
for each $1 \leq i \leq n$. Then there exists $N_{0} \in \N$ such that for all $N \geq N_{0}$ we have that
\begin{equation*}
\#\Qc(\bx, \bt, N) \geq \frac{1}{p}N^{n+1-\sum_{i=1}^{n}\tau_{i}}-1.
\end{equation*}
\end{lemma}
As with Theorem \ref{theorem1}, the equivalent version of this result in $\R^{n}$ has previously been proven, (see Lemma 3 of~\cite{RSS17}). Further, as $\sum_{i=1}^{n}\tau_{i}<n+1$ we can choose $N$ large enough such that
\begin{equation*}
\#\Qc(\bx,\bt, N) \geq \frac{1}{2p}N^{n+1-\sum_{i=1}^{n}\tau_{i}}.
\end{equation*}
Thus combining this with Theorem~\ref{theorem1} we have the expected result that $\#\Qc(\bx,\bt, N) \asymp N^{n+1-\sum_{i=1}^{n}\tau_{i}}$.\par 
The proofs of Lemma \ref{lemma1} and Lemma \ref{lemma2} use elementary techniques. The proof of Theorem \ref{theorem1} is more substantial and uses $p$-adic approximation lattices and lattice counting techniques. Prior to the proofs of these results we give an example of their applications in Diophantine approximation.
\section{$p$-adic Diophantine approximation}
 As an application of the main results in the previous section we consider the set of $p$-adic simultaneously approximable points over coordinate hyperplanes. Define the set of weighted simultaneously approximable points as follows. For an $n$-tuple of approximation functions $\Psi=(\psi_{1}, \dots, \psi_{n})$ and  $q_{0} \in \N$ let
 \begin{equation*}
 \Ap_{q_{0}}(\Psi)= \underset{1 \leq i \leq n}{\bigcup_{|q_{i}|\leq q_{0}, \, \, gcd(q_{i},q_{0})=1}}\left\{ \bx \in \Zp^{n}:\left|x_{i}-\frac{q_{i}}{q_{0}}\right|_{p}< \psi_{i}(q_{0}) \right\},
 \end{equation*}
 where $\bx=(x_{1}, \dots, x_{n})$. Define the set of weighted $\Psi$-approximable $p$-adic points as
 \begin{equation*}
 \Wp_{n}(\Psi):= \limsup_{q_{0} \to \infty} \Ap_{q_{0}}(\Psi).
 \end{equation*}
 If the approximation functions have the form $\Psi=(\psi, \dots , \psi)$ then we will use the notation $\Wp_{n}(\Psi)=\Wp_{n}(\psi)$, and if $\psi_{i}(q) = q^{-\tau_{i}}$ for each $1 \leq i \leq n$ and for some exponents of approximation $\bt=(\tau_{1}, \dots , \tau_{n}) \in \R^{n}_{>0}$ we will use the notation $\Wp_{n}(\Psi)=\Wp_{n}(\bt)$. By considering the Dirichlet style theorem for the set $\Wp_{n}(\bt)$ we have that $\Wp_{n}(\bt)= \Zp^{n}$ provided that $\sum_{i=1}^{n} \tau_{i} \leq n+1$. There are also a variety of results on the Haar measure of $p$-adic approximable points. The classical result was originally proven by Jarnik \cite{J45}, and since then there has been results in the linear forms \cite{L55}, weighted \cite{BLW20b}, and Duffin-Schaeffer \cite{H10} cases.
 \par
  For sets of zero Haar measure we use Hausdorff measure and Hausdorff dimension to provide more accurate notions of size. We briefly recap the definition and notation of Hausdorff measure and dimension. For a metric space $(X,d)$, a set $U \subset X$, and $\rho>0$, define a $\rho$-cover of $U$ as a sequence of balls $\{B_{i}\}$ such that $U \subset \bigcup_{i}B_{i}$ and for all balls $r(B_{i}) \leq \rho$, where $r(.)$ denotes the radius of the ball. Define a dimension function $f : \R_{+} \to \R_{+}$ as an increasing continuous function with $f(r) \to 0$ as $r \to 0$. Define the $f$-Hausdorff measure as
\begin{equation*}
\ha^{f}(U)= \lim_{\rho \to 0^{+}}\inf \left\{ \sum_{i}f(r(B_{i})): \{ B_{i} \} \, \, \text{ is a $\rho$-cover of} \, U \right\},
\end{equation*}
where the infimum is taken over all $\rho$-covers of $U$. When the dimension function $f(x)=x^{s}$ we will use the notation $\ha^{f}=\ha^{s}$. Define the Hausdorff dimension as
\begin{equation*}
\dim U = \inf \{ s\geq 0: \ha^{s}(U)= 0 \}.
\end{equation*}
 \par  
  In \cite{BLW20b} it was proven that, for $\sum_{i=1}^{n} \tau_{i} > n+1$,
 \begin{equation*}
 \dim \Wp_{n}(\bt)= \min_{1 \leq i \leq n} \left\{ \frac{n+1+ \sum_{j=i}^{n}(\tau_{i}-\tau_{j})}{\tau_{i}} \right\}.
 \end{equation*}
 It would be desirable to obtain equivalent measure results for simultaneous $p$-adic approximable points over manifolds. In \cite{KT07} Kleinbock and Tomanov proved the extremality of $p$-adic manifolds provided some non-degeneracy conditions are satisfied. Generally a manifold $\M \subset \Zp^{n}$ is said to be extremal if for almost all points $\bx \in \M$, with respect to the induced Haar measure of the manifold, we have that $\tau(\bx) = n+1$ (see \cite{KT07} for more details). There are a variety of results for $p$-adic dual approximation, see for example \cite{BBK05, BK03, B11, DG1, MG12}, however results in the simultaneous case are lacking. Recently Oliveira \cite{Oli20} produced a Khintchine-style Theorem for simultaneous $p$-adic approximation with denominators coming from $p$-adic balls. This result has a similar style to our result, with the difference being that our denominators come from a ball with radius tending to zero, rather than a fixed constant. Other than this there are relatively few Khintchine-style results. \par 
 For the Hausdorff dimension there are recent results on simultaneously approximable points over the Veronese curve, $\Cf=\{(x, x^{2}, \dots, x^{n}): x \in \Zp\}$, provided the approximation functions have sufficiently large Diophantine exponents \cite{BBS20, MR2837469}. In \cite{BLW20b} a lower bound for the Hausdorff dimension was found for general $n$-dimensional normal curves. A key reason the upper bound could not be obtained was a lack in results on the behaviour of rational points close to $p$-adic curves. The main results of this paper provide us with a good understanding of the behaviour of rational points close to coordinate hyperplanes. The results of this section are closely related to a variety of results in the real case on Diophantine approximation over coordinate hyperplanes, see \cite{BLVV17, R15, RSS17}. \par 
  For a $p$-adic integer $\bal \in \Zp^{m}$ for $1 \leq m \leq n-1$ define the coordinate hyperplane
 \begin{equation*}
 \Pi_{\bal}:= \{ (x_{1},\dots, x_{d}, \bal): (x_{1}, \dots , x_{d}) \in \Zp^{d} \} \subset \Zp^{n},
 \end{equation*}
 where $n=d+m$. For the set $\Wp_{n}(\bt) \cap \Pi_{\bal}$ we have the trivial result that
 \begin{equation*}
\dim \Wp_{n}(\bt) \cap \Pi_{\bal} \leq \dim \Pi_{\bal}=n-m,
\end{equation*} 
with equality when $\sum_{i=1}^{n} \tau_{i} \leq n+1$. In this paper we prove the following result on the Hausdorff dimension of $\W_{n}(\Psi) \cap \Pi_{\bal}$.
 \begin{theorem} \label{theorem2}
Let $\Pi_{\bal}$ be a coordinate hyperplane of $\Zp^{n}$, let $\bal \in \Zp^{m}$ satisfy $\tau(\bal)=m+1$. Let $\bt=(\tau_{1}, \dots, \tau_{n}) \in \R^{n}_{+}$ be a weight vector with the properties that
\begin{equation*}
\sum_{i=1}^{m}\tau_{d+i} < m+1, \quad \sum_{i=1}^{n} \tau_{i} > n+1, \quad \tau_{i}>1,
\end{equation*}
 for all $1 \leq i \leq n$. Then
\begin{equation*}
\dim \Wp_{n}(\bt) \cap \Pi_{\bal} = \min_{1 \leq i \leq d} \left\{ \frac{n+1-\sum_{i=1}^{m}\tau_{d+i}+\underset{1 \leq j \leq d}{\sum_{\tau_{j} \leq \tau_{i}}}(\tau_{i}-\tau_{j})}{\tau_{i}} \right\}=s.
\end{equation*}
Further
\begin{equation*}
\ha^{s}\left( \Wp_{n}(\bt) \cap \Pi_{\bal} \right)= \infty.
\end{equation*}
\end{theorem}
\begin{remark} \rm
The constraints on $(\tau_{d+1}, \dots, \tau_{n})$ ensure that we can apply Theorem \ref{theorem1}. The condition that $\sum_{i=1}^{n}\tau_{i}>n+1$ ensures that we do not include the trivial case when $\Wp_{n}(\bt)=\Zp^{n}$, in which case $\dim \Wp_{n}(\bt) \cap \Pi_{\bal}=n-m$.
\end{remark}
\begin{remark} \rm
 In the special case where the approximation functions are the same i.e. ($\bt=(\tau, \dots, \tau)$), then we have that, for $1+ \frac{1}{n} < \tau < 1+\frac{1}{m}$,
 \begin{equation*}
 \dim \W_{n}(\bt) \cap \Pi_{\bal}=\frac{n+1}{\tau}-m.
 \end{equation*}
 This gives us the expected dimension of the set of approximable points $\W_{n}(\bt)$ less the codimension of the hyperplane $\Pi_{\bal}$.
\end{remark}
\begin{remark}\rm
 We can use the same style of proof used to prove the upper bound of Theorem \ref{theorem2}, in combination with Lemma \ref{lemma1} rather than Theorem \ref{theorem1}, to prove that for any $\alpha \in \Zp$ and approximation exponent $\max\{1, \tau(\alpha)-1\}<\tau_{n} < \tau(\alpha)$ we have that
\begin{equation*}
\dim \Wp_{n}(\bt) \cap \Pi_{\alpha} \leq \min_{1 \leq i \leq n-1} \left\{ \frac{n+\tau(\alpha)-1-\tau_{n}+\underset{j \neq n}{\sum_{\tau_{j} \leq \tau_{i}}}(\tau_{i}-\tau_{j}),}{\tau_{i}} \right\}.
\end{equation*}
Proving the corresponding lower bound of this result is currently beyond our reach.
\end{remark}
 \par 
 For general approximation functions $\Psi=(\psi_{1}, \dots, \psi_{n})$, let
 \begin{equation} \label{limits}
 \psi_{i}^{*}= \lim_{ q \to \infty} \frac{-\log(\psi(q))}{\log q}.
 \end{equation}
 Providing the limits exists and are positive finite for each $1 \leq i \leq n$ then define $\Psi^{*}=(\psi^{*}_{1}, \dots, \psi_{n}^{*})$.
 \begin{corollary} \label{general5}
 Let $\Psi=(\psi_{1}, \dots, \psi_{n})$ be an $n$-tuple of approximation functions with each $\psi_{i}$ having positive finite limit \eqref{limits}. If $\Psi^{*}$ satisfy the same conditions as in Theorem \ref{theorem2}, then for all $\bal \in \Zp^{m}$ with $\tau(\bal)=m+1$,
 \begin{equation*}
 \dim \Wp_{n}(\Psi) \cap \Pi_{\bal} = {\min_{1 \leq i \leq d}} \left\{ \frac{n+1-\sum_{i=1}^{m}\psi^{*}_{d+i}+\underset{1 \leq j \leq d}{\sum_{\psi^{*}_{j} \leq \psi^{*}_{i}}}(\psi^{*}_{i}-\psi^{*}_{j})}{\psi^{*}_{i}} \right\}.
 \end{equation*}
 \end{corollary}
 The corollary easily follows from the observation that by the definition of \eqref{limits} there exists sufficiently large $q \in \N$ such that
 \begin{equation*}
 q^{-\psi^{*}_{i}-\epsilon_{i}} \leq \psi_{i}(q) \leq q^{-\psi^{*}_{i}+\epsilon_{i}}
 \end{equation*}
 for all $1 \leq i \leq n$ and $\epsilon=(\epsilon_{1}, \dots , \epsilon_{n})>0$ with $\epsilon_{i} \to 0$ as $q \to \infty$. And so 
 \begin{equation*}
 \Wp_{n}(\Psi^{*}+ \epsilon) \subseteq \Wp_{n}(\Psi) \subseteq \Wp_{n}(\Psi^{*}-\epsilon).
 \end{equation*}
 Letting $\epsilon \to 0$ we obtain the desired result. Note that while Corollary \ref{general5} provides a result for general $\Psi$ with components satisfying \eqref{limits}, there are many functions where such limits do not exist. \par
 
 The following section provides auxiliary results needed to prove Theorem \ref{theorem2}. In particular the framework for the Mass Transference Principle from rectangles to rectangles \cite{WW19} is provided. This result is crucial in finding the lower bound result of Theorem \ref{theorem2}.
\section{Auxiliary results} \label{aux}
We provide a brief set of known results that we will use in the proof of Theorem \ref{theorem2}. The first result we state can be considered as the $p$-adic version of Minkowski's theorem for systems of linear forms. The proof is a straightforward application of the pigeon-hole principle and can be found in \cite{BLW20b}.
\begin{lemma} \label{mink5}
Let $L_{i}(\bx):\Zp^{n} \to \Zp$, with $i=1, \dots , n$, be linear forms with $p$-adic integer coefficients.  Let $\sum_{i=1}^{n} \tau_{i}=n+1$ for $\tau_{i} \in \R_{+}$, and $H \geq 1$. Then there exists a non-zero rational integer vector $\bx=(x_{0},x_{1}, \dots , x_{n})$ with
\begin{equation*}
\max_{0 \leq i \leq n} |x_{i}| \leq H 
\end{equation*}
satisfying the system of inequalities
\begin{equation*}
|L_{i}(\bx)|_{p} < pH^{-\tau_{i}} \, \, \text{ for } i=1, \dots , n.
\end{equation*}
\end{lemma}
The following lemma generally states that the measure of a $\limsup$ set of balls remains unaltered when the radius is multiplied by some constant. The Euclidean version of this result is well known and appears in a variety of texts, see \cite{BDV06}. The following version for ultrametric spaces was proven in \cite{BLW20b}.
\begin{lemma}\label{measure_unchanged1}
Let $(X,d)$ be a separable ultrametric space and $\mu$ be a Borel regular measure on $X$. Let $(B_i)_{i\in\N}$ be a sequence of balls in $X$ with radii $r_i\to0$ as
$i\to\infty$. Let $(U_i)_{i\in\N}$ be a sequence of $\mu$-measurable
sets such that $U_i\subset B_i$ for all $i$. Assume that for some
$c>0$
\begin{equation*}
     |U_i|\ge c|B_i|\qquad\text{for all }i\,.
\end{equation*}
Then the limsup sets
$$
\textstyle
 \mathcal{U}=\limsup\limits_{i\to\infty}U_i:=\bigcap\limits_{j=1}^\infty\ \bigcup\limits_{i\ge
 j}U_i\qquad\text{ and }\qquad
 \mathcal{B}=\limsup\limits_{i\to\infty}B_i:=\bigcap\limits_{j=1}^\infty\ \bigcup\limits_{i\ge
 j}B_i
$$
have the same $\mu$-measure.
\end{lemma}
In particular, if we chose the approximation function $\psi_{i}(q)=pq^{1+1/n}$ for each $1 \leq i \leq n$ then by Lemma \ref{mink5} we know $\Wp_{n}(\Psi)=\Zp$, and so by shrinking the $\limsup$ set of balls by constant $1/p$ Lemma \ref{measure_unchanged1} gives us that $\mu_{n}(\Wp_{n}(\Psi/p))=1$. \par 
 
 Another key result in our proof of Theorem \ref{theorem2} is the following Mass Transference Principle type theorem. In order to state this theorem we need the notion of local ubiquity for rectangles, a variation of the notion of ubiquity introduced by Beresnevich, Dickinson, and Velani \cite{BDV06}. Fix an integer $n \geq 1$, and for each $1 \leq i \leq n$ let $(X_{i},|.|_{i}, m_{i})$ be a bounded locally compact metric space with $m_{i}$ a $\delta_{i}$-Ahlfors probability measure. Consider the product space $(X, |.|,m)$, where
\begin{equation*}
X=\prod_{i=1}^{n}X_{i}, \quad m=\prod_{i=1}^{n}m_{i}, \quad |.|=\max_{1 \leq i \leq n}|.|_{i}.
\end{equation*}
For any $x \in X$ and $r \in \R_{+}$ define the open ball
\begin{equation*}
B(x,r)=\left\{ y \in X: \max_{1 \leq i \leq n}|x_{i}-y_{i}|_{i}< r \right\}=\prod_{i=1}^{n}B_{i}(x_{i},r),
\end{equation*}
where $B_{i}$ are the usual balls associated with the $i^{\text{th}}$ metric space. Let $J$ be a infinite countable index set, and $\beta: J \to \R_{+}$ a positive function. Let $l_{k},u_{k}$ be two sequences in $\R_{+}$ such that $u_{k} \geq l_{k}$ with $l_{k} \to \infty$ as $k \to \infty$. Define
\begin{equation*}
J_{k}= \left\{ \alpha \in J: l_{k} \leq \beta_{\alpha} \leq u_{k} \right\}.
\end{equation*}
Let $\rho: \R_{+} \to \R_{+}$ be a non-increasing function with $\rho(\beta_{\alpha}) \to 0$ as $\beta_{\alpha} \to \infty$. For each $1 \leq i \leq n$, let $\{ R_{\alpha,i}\}_{\alpha \in J}$ be a sequence of subsets in $X_{i}$. As with the standard setting of ubiquitous systems define the resonant sets
\begin{equation*}
\left\{ R_{\alpha}=\prod_{i=1}^{n} R_{\alpha, i} \right\}_{ \alpha \in J}.
\end{equation*}
For $\ba=(a_{1}, \dots, a_{n}) \in \R_{+}^{n}$ denote the set of hyperrectangles
\begin{equation*}
\Delta(R_{\alpha},\rho(r)^{\ba})= \prod_{i=1}^{n}  \Delta(R_{\alpha,i},\rho(r)^{a_{i}}),
\end{equation*}
where for some set $A$ and $b \in \R_{+}$
\begin{equation*}
\Delta(A,b)= \bigcup_{a \in A}B(a,b).
\end{equation*}
\begin{definition}[local ubiquitous system of rectangles]
Call $(\{R_{\alpha}\}_{\alpha \in J}, \beta)$ a local ubiquitous system of rectangles with respect to $(\rho,\ba)$ if there exists a constant $c>0$ such that for any ball $B \subset X$,
\begin{equation*}
\limsup_{k \to \infty} m \left( B \cap \bigcup_{\alpha \in J_{k}}\Delta(R_{\alpha}, \rho(u_{k})^{\ba}) \right) \geq c m(B).
\end{equation*}
\end{definition}
The second property needed to state the theorem is the following local scaling property, first introduced in \cite{AB19}. While we will not need it in our use of the Theorem \ref{MTPRR5} we state it and include it in the final theorem for completeness.
\begin{definition}[$k$-scaling property] \label{LSP}
Let $0 \leq k <1$ and $1 \leq i \leq n$. Then $\{R_{\alpha,i}\}_{\alpha \in J}$ has $k$-scaling property if for any $\alpha \in J$, any ball $B(x_{i},r) \subset X_{i}$ with centre $x_{i} \in R_{\alpha,i}$, and $0 < \epsilon < r$ then
\begin{equation*}
c_{2}r^{\delta_{i}k}\epsilon^{\delta_{i}(1-k)} \leq m_{i} \left( B(x_{i},r) \cap \Delta(R_{\alpha,i},\epsilon) \right) \leq c_{3} r^{\delta_{i}k}\epsilon^{\delta_{i}(1-k)},
\end{equation*}
for some constants $c_{2},c_{3}>0$.
\end{definition}
In our use the resonant sets will be sets of points, so $k=0$. For $\textbf{t}=(t_{1}, \dots, t_{n}) \in \R^{n}_{+}$ define
\begin{equation*}
W(\textbf{t})= \limsup_{ \alpha \in J} \Delta \left( R_{\alpha},\rho(\beta_{\alpha})^{\ba+\textbf{t}} \right).
\end{equation*}
Given the above notations and definitions we can state the Mass Transference Principle from rectangles to rectangles (MTPRR) of \cite{WW19}.
\begin{theorem}[Mass Transference Principle from rectangles to rectangles] \label{MTPRR5}
Under the settings above assume that $(\{ R_{\alpha} \}_{\alpha \in J}, \beta)$ satisfies the local ubiquity for systems of rectangles condition with respect to $(\rho, \ba)$, and the $k$-scaling property. Then
\begin{equation*}
\dim W(\textbf{t}) \geq \min_{A_{i} \in A} \left\{ \sum_{j \in K_{1}} \delta_{j}+ \sum_{j \in K_{2}} \delta_{j}+ k \sum_{j \in K_{3}} \delta_{j}+(1-k) \frac{\sum_{j \in K_{3}}a_{j}\delta_{j}-\sum_{j \in K_{2}}t_{j}\delta_{j}}{A_{i}} \right\}=s,
\end{equation*}
where $A=\{ a_{i}, a_{i}+t_{i} , 1 \leq i \leq n \}$ and $K_{1},K_{2},K_{3}$ are a partition of $\{1, \dots, n\}$ defined as
\begin{equation*}
 K_{1}=\{ j:a_{j} \geq A_{i}\}, \quad K_{2}=\{j: a_{j}+t_{j} \leq A_{i} \} \backslash K_{1}, \quad K_{3}=\{1, \dots n\} \backslash (K_{1} \cup K_{2}).
 \end{equation*}
 Further,
 \begin{equation*}
 \ha^{s}(B \cap W(\bt))=\ha^{s}(B).
 \end{equation*}
 \end{theorem}
Hence, provided we can find a $\limsup$ set of hyperrectangles that satisfy the local ubiquity property for rectangles then we have a lower bound for the corresponding $\limsup$ set of shrunken hyperrectangles.

\section{Proof of Theorem \ref{theorem2}}
We split the proof into the upper and lower bound, and solve each case separately. In both cases we will use the following simplified set. Let $\pi$ be the projection $\pi: \Zp^{n} \to \Zp^{n-m}$, defined by
\begin{equation*}
(x_{1}, \dots , x_{n}) \mapsto (x_{1}, \dots , x_{d}).
\end{equation*}
 By a well known theorem of Hausdorff theory (see Proposition 3.3 of \cite{F14}) as $\pi$ is a bi-Lipschitz mapping over $\Wp_{n}(\bt)\cap \Pi_{\bal}$, we have that
\begin{equation*}
\dim \Wp_{n}(\bt) \cap \Pi_{\bal} = \dim \pi(\Wp_{n}(\bt) \cap \Pi_{\bal}).
\end{equation*}
Let $\bt_{m}=(\tau_{d+1}, \dots , \tau_{n})$ denote the $m$-tuple of approximation exponents over $\bal$ and similarly let $\bt_{d}=(\tau_{1}, \dots , \tau_{d})$ denote the $d$-tuple of approximation exponents over the independent variables of $\Pi_{\bal}$. Consider the set of integers
\begin{equation*}
\Qc(\bal, \bt_{m}):= \left\{ q_{0} \in \N : \left|\alpha_{i} - \frac{q_{d+i}}{q_{0}} \right|_{p} < q_{0}^{-\tau_{d+i}}, \, \text{ for some } \, \begin{array}{c} |q_{i}| \leq q_{0}, \\
gcd(q_{i},q_{0})=1,\\
\end{array} \,  1 \leq i \leq m \right\},
\end{equation*}
and the union of sets
\begin{equation*}
\Ap_{q_{0}}^{*}(\bt_{d})= \underset{1 \leq i \leq d}{\bigcup_{|q_{i}| \leq q_{0}, \, \, gcd(q_{i},q_{0})=1}} \left\{ \bx \in \Zp^{d}: \left|x_{i}-\frac{q_{i}}{q_{0}} \right|_{p} < q_{0}^{-\tau_{i}} \right\}.
\end{equation*}
Then,
\begin{equation*}
\pi(\W_{n}(\bt) \cap \Pi_{\bal})=\limsup_{q_{0} \in \Qc(\bal, \bt_{m})}\Ap^{*}_{q_{0}}(\bt_{d}),
\end{equation*}
 hence we only need to find the upper and lower bounds for $\dim \underset{q_{0} \in \Qc(\bal, \bt_{m})}{\limsup} \A^{*}_{q_{0}}(\bt_{d})$.

 \subsection{Upper bound}
For the upper bound we take the standard cover of hyperrectangles used in the construction of $\Ap^{*}_{q}(\bt_{d})$. By a standard geometrical argument note that each hyperrectangle, centred at some $\left(\frac{q_{1}}{q}, \dots , \frac{q_{n}}{q_{0}} \right) \in \Q^{d}$ in the construction of $\Ap_{q}^{*}(\bt_{d})$, can be covered by a finite collection of balls $\mathfrak{B}_{q}(\tau_{i})$ of radius $q^{-\tau_{i}}$ for $1 \leq i \leq d$. Without loss of generality we can assume that
\begin{equation*}
\tau_{1} \geq \dots \geq \tau_{d},
\end{equation*}
since if not then we could take some bi-Lipschitz mapping to reorder the coordinate axes such that this was the case. Hence for each $j \leq i$,
\begin{equation*}
\frac{q^{-\tau_{j}}}{q^{-\tau_{i}}} \leq 1.
\end{equation*}
Hence in the product below we only consider the $j \geq i$. By the above argument we have that the cardinality of $\mathfrak{B}_{q}(\tau_{i})$ is
\begin{equation*}
\# \mathfrak{B}_{q}(\tau_{i}) \ll \prod_{j=i}^{d} \frac{q^{-\tau_{j}}}{q^{-\tau_{i}}}=q^{\sum_{j=i}^{d}(\tau_{i}-\tau_{j})}.
\end{equation*}
 As each $\tau_{i}$-approximation function is decreasing as $q$ increases, for each interval $2^{k} \leq q < 2^{k+1}$ take $q=2^{k}$ over such interval. Let 
 \begin{equation*}
\Qc^{'}(\bx, \bt_{m}, N):= \left\{ q_{0}\in \N: (q_{0}, \dots , q_{m}) \in \Qc(\bx, \bt_{m}, N) \, \text{ and } gcd(q_{i},q_{0})=1 \right\}.
\end{equation*}
Since each $\tau_{i}>1$ for $1 \leq i \leq m$ each $q_{0}$ has unique associated $(q_{1}, \dots, q_{m})$ in $\Qc(\bx, \bt_{m},N)$ so we have that $\#\Qc^{'}(\bx, \bt_{m},N) \leq \# \Qc(\bx, \bt_{m}, N)$. Further, by the coprimality of each $q_{i}$ with $q_{0}$ note that the inequalities
\begin{equation*}
|q_{0}x_{i}-q_{i}|_{p}<H^{-\tau_{i}}, \quad \text{and} \quad \left|x_{i}-\frac{q_{i}}{q_{0}} \right|_{p}<H^{-\tau_{i}}
\end{equation*}
are equivalent since $p \not | q_{0}$. To check this observe that each $x_{i} \in \Zp$ and then use the strong triangle inequality. \par 
 Given the above we have that $\Qc(\bal , \bt_{m}) \subseteq \bigcup_{k \in \N} \Qc^{'}(\bal, \bt_{m}, 2^{k})$. Hence
\begin{align*}
\ha^{s}\left( \limsup_{q \in \Qc(\bal, \bt_{m})}\A^{*}_{q}(\bt_{d})\right) & \leq \sum_{k=1}^{\infty} \sum_{q \in \Qc(\bal, \bt_{m}, 2^{k})} \phi(q)^{d} \#\mathfrak{B}_{q}(\tau_{i}).(q^{-\tau_{i}})^{s}, \\
& \overset{(a)}{\ll} \sum_{k=1}^{\infty} 2^{k(m+1-\sum_{i=1}^{m}\tau_{d+i})}(2^{k+1})^{d}(2^{k+1})^{\sum_{j=i}^{d}(\tau_{i}-\tau_{j})}(2^{k})^{-\tau_{i}s}, \\
& \ll \sum_{k=1}^{\infty}2^{k(n+1-\sum_{i=1}^{m}\tau_{d+i}+\sum_{j=i}^{d}(\tau_{i}-\tau_{j})-\tau_{i}s)},
\end{align*}
where $(a)$ follows from Theorem \ref{theorem1}. The above sum converges when
\begin{equation*}
s\geq \frac{n+1-\sum_{i=1}^{m}\tau_{d+i}+\sum_{j=i}^{d}(\tau_{i}-\tau_{j})}{\tau_{i}} + \epsilon,
\end{equation*}
for any $\epsilon>0$. This is true for each $1 \leq i \leq d$, and as $\epsilon$ is arbitrary, we have that 
\begin{equation*}
s\geq \min_{1 \leq i \leq d} \left\{ \frac{n+1-\sum_{i=1}^{m}\tau_{d+i}+\sum_{j=i}^{d}(\tau_{i}-\tau_{j})}{\tau_{i}} \right\},
\end{equation*}
completing the upper bound result. Note that the result of Remark 2.2 can similarly be obtained by replacing Theorem \ref{theorem1} used at $(a)$ by Lemma \ref{lemma1}.
\subsection{Lower bound}
In order to use Theorem \ref{MTPRR5} to prove the lower bound of Theorem \ref{theorem2} we need to construct a ubiquitous system of rectangles. In following with the ubiquity setup for Theorem \ref{MTPRR5} let
\begin{equation*}
\begin{array}{ccc}
J=\Qc(\bal, \bt_{m}), & R_{q,i}=\left\{ \frac{q_{i}}{q} \in \Q : \begin{array}{c} |q_{i}| \leq q, \\
gcd(q_{i},q)=1
\end{array} \right\}, & R_{q}= \prod_{i=1}^{d}R_{q,i},\\
\beta(q)=q, & \rho(q)=q^{-1}, & l_{k}=M^{k}, \, \, u_{k}=M^{k+1},
\end{array}
\end{equation*} 
where $M \in \N$ is a fixed integer to be determined later. Then we have that
\begin{equation*}
J_{k}= \{ q \in \Qc(\bal, \bt_{m}): M^{k} \leq q < M^{k+1} \}.
\end{equation*}
Note that $J_{k} \subseteq \Qc^{'}(\bal,\bt_{m},2^{k+1})$.
For a vector $\ba=(a_{1}, \dots , a_{n}) \in \R^{n}_{+}$ let
\begin{equation*}
\Delta(R_{q}, \rho(r)^{\ba})= \prod_{i=1}^{n} \bigcup_{q_{i} \in R_{q,i}}B\left( \frac{q_{i}}{q}, r^{-a_{i}} \right).
\end{equation*}
We prove the following.
\begin{proposition} \label{ubiquity}
Let $R_{q}$, $\rho$, and $J_{k}$ be as above, and let $\tv=(v_{1}, \dots, v_{d}) \in \R^{d}_{>0}$ with each $v_{i} > 1$ and  
\begin{equation*}
\sum_{i=1}^{d}v_{i}=n+1-\sum_{i=1}^{m}\tau_{d+i},
\end{equation*}
for $\sum_{i=1}^{m}\tau_{d+i} < m+1$ and each $\tau_{i}>1$. Then for any ball $B=B(x,r) \subset \Zp^{d}$, with centre $x \in \Zp^{d}$ and radius $0<r<r_{0}$ for some $r_{0} \in \R_{+}$, there exists a constant $c>0$ such that
\begin{equation*}
\mu_{p,d}\left(B \cap \bigcup_{q \in J_{k}} \Delta(R_{\bq}, \rho(u_{k})^{\tv})  \right) \geq c \mu_{p,d}(B),
\end{equation*}
provided $M > (3^{d}C_{1})^{\frac{1}{n+1-\sum_{i=1}^{m}\tau_{d+i}}}$.
\end{proposition}
The proof of this result follows the same style of many similar results in $\R^{n}$. For example see Theorem 1.3 of \cite{BRV16} for the one dimensional real case, or Proposition 5.1 of \cite{BLW20b} for the $n$-dimensional $p$-adic case.
\begin{proof} 
For any $y=(y_{1}, \dots, y_{d}) \in (\Zp \backslash \Q)^{d}$, consider the system of inequalities
\begin{equation} \label{system}
\begin{cases}
|q_{0}\alpha_{i}-q_{d+i}|_{p} < (M^{k+1})^{-\tau_{d+i}}, \quad 1 \leq i \leq m, \\
|q_{0}y_{i}-q_{i}|_{p}< (M^{k+1})^{-v_{i}}, \quad 1 \leq i \leq d, \\
\max_{0 \leq i \leq n} |q_{i}| \leq M^{k+1}.
\end{cases}
\end{equation}
By the condition on $\tv$ we have, by Lemma \ref{mink5}, that there exists a non-zero integer solution $(q_{0}, \dots , q_{n}) \in \Z^{n+1}$ to \eqref{system}. Assume without loss of generality that $q_{0} \geq 0$. We prove that there exists a rational integer solution $(q_{0}, \dots, q_{n})$ to \eqref{system} satisfying
\begin{enumerate}[i)]
\item $q_{0} \neq 0$, 
\item $\frac{q_{i}}{q_{0}} \in \Zp$ for each $1 \leq i \leq n$, 
\item $gcd(q_{0},p)=1$.
\end{enumerate}
Firstly, suppose $q_{0}=0$, then the above equations imply that each $|q_{i}|_{p}<(M^{k+1})^{-\tau_{d+i}}$ for $1 \leq i \leq m$ or $|q_{i}|_{p}<(M^{k+1})^{-v_{i}}$ for $1 \leq i \leq d$. However, since each $\tau_{i}$ and $v_{i}$ are greater than one such inequalities are impossible, unless $q_{i}=0$ for all $1 \leq i \leq n$. But this solution is identically zero. \par 
For ii) note that since $q_{0} \neq 0$ we may divide each inequality in \eqref{system} by $q_{0}$ to get
\begin{equation*}
\begin{cases}
|q_{0}|_{p}\left|\alpha_{i}-\frac{q_{d+i}}{q_{0}}\right|_{p} < (M^{k+1})^{-\tau_{d+i}}, \quad 1 \leq i \leq m, \\
|q_{0}|_{p} \left|y_{i}-\frac{q_{i}}{q_{0}} \right|_{p}< (M^{k+1})^{-v_{i}}, \quad 1 \leq i \leq d, \\
\max_{0 \leq i \leq n} |q_{i}| \leq M^{k+1}.
\end{cases}
\end{equation*}
Suppose $\left| \frac{q_{i}}{q_{0}}\right|_{p}>1$. Noting that each $y_{i}$ and $\alpha_{i}$ are $p$-adic integers, and using the strong triangle inequality, we obtain that each $|q_{i}|_{p}<(M^{k+1})^{-\tau_{i}}$ for $1 \leq i \leq m$ and $|q_{i}|_{p}<(M^{k+1})^{-v_{i}}$ for $1 \leq i \leq d$. As previously stated such condition is impossible unless each $q_{i}=0$, in which case $\left| \frac{q_{i}}{q_{0}}\right|_{p}=0$, contradicting that $\left| \frac{q_{i}}{q_{0}}\right|_{p}> 1$. \par 
For iii) suppose that $(q_{0}, \dots , q_{n})$ is a solution to \eqref{system} and suppose that $p^{t}|q_{0}$ but $p^{t+1} \not | q_{0}$. Then by ii) we have that $p^{t}|q_{i}$ for all $1 \leq i \leq n$. Let $q_{i}'=p^{-t}q_{i}$ for each $0 \leq i \leq n$. Note that 
\begin{equation*}
\max_{0 \leq i \leq n} |q_{i}'| \leq \frac{M^{k+1}}{p^{k}}=H^{'}.
\end{equation*}
Then for each $1 \leq i \leq m$ we have that
\begin{align*}
|q_{0}'\alpha_{i}-q_{i}'|_{p} & = p^{k}p^{-k}|q_{0}'\alpha_{i}-q_{i}'|_{p}, \\
&=p^{k}|q_{0}\alpha_{i}-q_{i}|_{p}, \\
&<p^{k}(M^{k+1})^{-\tau_{i}}, \\
&<(H^{'})^{-\tau_{i}},
\end{align*}
 and similarly for the approximations over $y$. Hence $(q_{0}', \dots , q_{n}')$ is a rational integer solution with $gcd(q_{0}',p)=1$. Henceforth we will suppose $p \not |q_{0}$. \par 
 Since $gcd(q_{0},p)=1$ we may divide \eqref{system} through by $|q_{0}|_{p}=1$ and the set of possible $y \in \Zp^{d}$ remain unchanged. Lastly, note that if $q_{0}$ has an associated solution to \eqref{system}, then $q_{0} \in \Qc^{'}(\bal,\bt_{m}, M^{k+1})$, thus we have that
\begin{equation*}
\mu_{p,d} \left( B \cap \bigcup_{q_{0} \in \Qc^{'}(\bal,\bt_{m},M^{k+1})} \Delta( R_{q_{0}}, \rho(M^{k+1})^{\tv}) ) \right)= \mu_{p,d}(B)
\end{equation*}
for any $B \subseteq \Zp^{d}$. Note that
\begin{align*}
\mu_{p,d} \left( B \cap \bigcup_{q_{0} \in \Qc^{'}(\bal,\bt_{m},M^{k+1})} \Delta( R_{q_{0}}, \rho(M^{k+1})^{\tv}) ) \right) \leq \mu_{p,d} & \left(  B \cap \bigcup_{q_{0} \in \Qc^{'}(\bal, \bt_{m}, M^{k})} \Delta( R_{q_{0}}, \rho(M^{k+1})^{\tv}) ) \right) \\
& + \mu_{p,d} \left( B \cap \bigcup_{q_{0} \in J_{k}} \Delta( R_{q_{0}}, \rho(M^{k+1})^{\tv}) ) \right),
\end{align*}
and so
\begin{equation*}
\mu_{p,d} \left( B \cap \bigcup_{q_{0} \in J_{k}} \Delta( R_{q_{0}}, \rho(M^{k+1})^{\tv}) ) \right) \geq \mu_{p,d}(B)- \mu_{p,d} \left( B \cap \bigcup_{q_{0} \in \Qc^{'}(\bal, \bt_{m}, M^{k})} \Delta( R_{q_{0}}, \rho(M^{k+1})^{\tv}) ) \right).
\end{equation*}
At this point we only want the $\frac{\bq}{q_{0}}=\left( \frac{q_{1}}{q_{0}}, \dots , \frac{q_{d}}{q_{0}} \right) \in R_{q_{0}}$ such that 
\begin{equation*}
B \cap B \left( \frac{\bq}{q}, \rho(M^{k+1})^{\tv} \right) \neq \emptyset.
\end{equation*}
For ball $B=B(x,r)$ with $x \in \Zp^{d}$ and $r \in \{p^{j}:j \in \Z \}$, this is equivalent to the set of solutions to
\begin{equation} \label{system2}
\left| x_{i}-\frac{q_{i}}{q_{0}} \right|_{p}<r, \quad 1 \leq i \leq d.
\end{equation}
For $q_{0}$ fixed and each $|q_{i}| \leq q_{0}$ by congruence classes we have that there are at most 
\begin{equation*}
(2q_{0}r+1)^{d}
\end{equation*}
suitable values of $\bq$. We can choose suitably large $k \in \N$ such that $M^{k}r>1$, and so for each $|q_{i}| \leq M^{k}$, $1 \leq i \leq d$ there are at most
\begin{equation} \label{no_q}
(3M^{k}r)^{d}
\end{equation}
possible values of $\bq$ solving \eqref{system2}.
 Hence
\begin{align*}
\mu_{p,d} \left( B \cap \bigcup_{q_{0} \in \Qc^{'}(\bal, \bt_{m}, M^{k})}\Delta( R_{q_{0}}, \rho(M^{k+1})^{\tv}) ) \right) & \leq \sum_{q_{0} \in \Qc^{'}(\bal, \bt_{m}, M^{k})} \sum_{\bq \, \, \text{solving} \, \, \eqref{system2}} \mu_{p,d} \left( B \cap \Delta\left( \frac{\bq}{q_{0}}, \rho(M^{k+1})^{\tv}) \right) \right), \\
& \overset{(a)}{\leq} \sum_{q_{0} \in \Qc^{'}(\bal, \bt_{m}, M^{k})} (3M^{k}r)^{d}(M^{k+1})^{-d\sum_{i=1}^{d}v_{i}}, \\
& \overset{(b)}{\leq} C_{1}M^{k(m+1-\sum_{i=1}^{m}\tau_{d+i})}3^{d}M^{kd}M^{-(k+1)(n+1-\sum_{i=1}^{m}\tau_{d+i}}\mu_{p,d}(B),\\
& \leq 3^{d}C_{1} M^{-n-1+\sum_{i=1}^{m}\tau_{d+i}}\mu_{p,d}(B),
\end{align*}
where $(a)$ follows by \eqref{no_q} and $(b)$ follows by Theorem \ref{theorem1} and our condition on $\tv$. As $M> (3^{d}C_{1})^{\frac{1}{n+1-\sum_{i=1}^{m}\tau_{d+i}}} $,
\begin{equation*}
c=\left( 1- \frac{3^{d}C_{1}}{M^{n+1-\sum_{i=1}^{m}\tau_{d+i}}} \right)>0.
\end{equation*}
Thus,
\begin{equation*}
\mu_{p,d} \left( B \cap \bigcup_{q_{0} \in J_{k}} \Delta( R_{q_{0}}, \rho(M^{k+1})^{\tv}) ) \right) \geq  c.\mu_{d}(B).
\end{equation*}
\end{proof}
Given Proposition \ref{ubiquity} we have that $(R_{q}, \beta)$ is a local ubiquitous system of rectangles with respect to $(\rho,\tv)$, provided $\sum_{i=1}^{d}v_{i}= n+1-\sum_{i=1}^{m}\tau_{d+i}$. Given $\bt_{d}=(\tau_{1}, \dots, \tau_{d}) \in \R^{d}_{>0}$ assume without loss of generality that $\tau_{1}> \tau_{2}> \dots > \tau_{d}$ and define each $v_{d-i}$ recursively by
\begin{equation*}
v_{d-i}=\min \left\{ \tau_{d-i}, \frac{n+1-\sum_{i=1}^{m}\tau_{d+i}-\sum_{j=d-i+1}^{d}v_{j}}{d-i} \right\}.
\end{equation*}
By the condition on $\bt_{d}$ of Theorem \ref{theorem2}, there exists a $k \in \{1, \dots, d\}$ such that
\begin{equation*}
v_{l}=\frac{n+1-\sum_{i=1}^{m}\tau_{d+i}-\sum_{j=d-k+1}^{d}v_{j}}{d-k},
\end{equation*}
for all $1 \leq l \leq d-k$. Clearly each $v_{i}\leq\tau_{i}$ for $1 \leq i \leq d$, and so 
the associated vector $\mathbf{t}=(t_{1}, \dots t_{n-1}) \in \R^{n-1}_{\geq 0}$ is defined by
\begin{equation*}
t_{i}=\tau_{i}-v_{i}, \quad 1 \leq i \leq d.
\end{equation*}
Consider the set
\begin{equation*}
A=\{v_{1}, \dots , v_{d}, \tau_{1}, \dots , \tau_{d} \}.
\end{equation*}
For each $A_{i} \in A$ observe the following:
\begin{enumerate}[i)]
\item $A_{i} \in \{v_{1}, \dots , v_{d}\}$: Then we have the sets
\begin{equation*}
\begin{array}{ccc}
K_{1}=\{1, \dots , \max\{i,d-k\} \}, & K_{2}=\{ \max\{i+1,d-k+1\}, \dots , d\}, & K_{3}=\emptyset.
\end{array}
\end{equation*}
By Theorem \ref{MTPRR5} we have that
\begin{align*}
\dim \Wp_{n}(\bt)\cap \Pi_{\bal} & \geq \min_{A_{i}} \left\{ \frac{\max\{i,d-k\}v_{i}+(d-\max\{i+1,d-k+1\})v_{i}- \sum_{j=\max\{i+1,d-k+1\}}^{d}t_{j}}{v_{i}} \right\}, \\
& =  \min_{A_{i}} \left\{ \frac{d v_{i}- \sum_{j=\max\{i+1,d-k+1\}}^{d}t_{j}}{v_{i}} \right\}.
\end{align*}
Since $t_{j}=0$ for $d-k+1 \leq j \leq d$ the above equation gives that $\dim \Wp_{n}(\bt)=d=n-m$, the maximal dimension of $\Wp_{n}(\bt) \cap \Pi_{\bal}$. 
\item $A_{i} \in \{\tau_{1}, \dots , \tau_{d}\}$: Since $\tau_{i}=v_{i}$ for $d-k+1 \leq i \leq d$ the above argument covers such case, so we only need to consider $\tau_{i}$ for $1 \leq i \leq d-k$. For such $\tau_{i}$ we have the sets
\begin{equation*}
\begin{array}{ccc}
K_{1}=\emptyset, & K_{2}=\{ i, \dots, d\}, & K_{3}=\{1, \dots , i-1\}.
\end{array}
\end{equation*}
Applying Theorem \ref{MTPRR5} we have
\begin{align*}
\dim \Wp_{n}(\bt) \cap \Pi_{\bal} & \geq \min_{A_{i}} \left\{ \frac{(d-i)\tau_{i} + \sum_{j=1}^{i-1}v_{j} - \sum_{j=i}^{d}t_{j}}{\tau_{i}} \right\}, \\
& =  \min_{A_{i}} \left\{ \frac{(d-i)\tau_{i}+ (d-k)\left( \frac{n+1-\sum_{i=1}^{m}\tau_{d+i}-\sum_{j=d-k+1}^{d}v_{j}}{d-k} \right)-\sum_{j=1}^{d-k}\tau_{j}}{\tau_{i}} \right\}, \\
&= \min_{A_{i}} \left\{ \frac{n+1-\sum_{i=1}^{m}\tau_{d+i}+\sum_{j=i}^{d}(\tau_{i}-\tau_{j})}{\tau_{i}} \right\}.
\end{align*}
\end{enumerate}
Combining i) and ii) we have that
\begin{equation*}
\dim \Wp_{n}(\bt) \cap \Pi_{\bal} \geq \min_{1 \leq i \leq d} \left\{ \frac{n+1-\sum_{i=1}^{m}\tau_{d+i}+\sum_{j=i}^{d}(\tau_{i}-\tau_{j})}{\tau_{i}} \right\},
\end{equation*}
completing the proof.

\section{Proof of the counting results}
Recall, we are aim to provide bounds on the set
\begin{equation*}
\Qc(\bx,\Psi, N):=\left\{ (q_{0}, \dots, q_{n}) \in \Z^{n+1} : \begin{array}{c}
0<q_{0} \leq N, \\
\max_{1 \leq i \leq n}|q_{i}| \leq N, \\
\end{array}
    \, \left| q_{0}x_{i}- q_{i} \right|_{p} < \psi_{i}(N), \, 1 \leq i \leq n \right\}.
\end{equation*}
We begin with the proof of Lemma \ref{lemma2}. This style of proof is not new and follows a similar method to the proof in the euclidean case (see Lemma 3 of \cite{RSS17}). \\
\vspace{1ex}

\textit{Proof of Lemma \ref{lemma2}:}
Fix $\bx=(x_{1}, \dots , x_{n}) \in \Zp^{n}$ and take $t=(t_{1}, \dots , t_{n}) \in \N^{n}$ to be the integers such that
\begin{equation*}
p^{-t_{i}} \leq N^{-\tau_{i}}<p^{-t_{i}+1}, \quad 1 \leq i \leq n.
\end{equation*}
Denote by $P=\prod_{i=1}^{n}p^{t_{i}}$. Consider a set of open disjoint rectangles $\{R_{i} \}_{i=1}^{P}$, each with some centre point $k_{i}=(k_{i,1}, \dots, k_{i,n}) \in \Z^{n}$ and sidelenghts $p^{-t_{i}}$. Choose the set of points $\{k_{i}\}$ such that $\Zp^{n} \subseteq \bigcup_{i=1}^{P}R_{i}$. Consider the $(N+1)^{n+1}$ set of points of the form
\begin{equation*}
(q_{0}x-q)=(q_{0}x_{1}-q_{1}, \dots , q_{0}x_{n}-q_{n}) \in \Zp^{n},
\end{equation*}
with $q_{i} \in [0,N]$ for each $0 \leq i \leq n$. By the Pigeon-hole principle there exists at least one rectangle, say $R_{j}$, containing at least 
\begin{equation*}
\frac{(N+1)^{n+1}}{P} > \frac{1}{p^{n}}N^{n+1-\sum_{i=1}^{n}\tau_{i}}
\end{equation*}
points. As $\sum_{i=1}^{n}\tau_{i}<n+1$ we can choose $N$ sufficiently large enough such that $p^{-n}N^{n+1-\sum_{i=1}^{n}\tau_{i}}>2$. Order the points $(q_{0},\dots ,q_{n})$, correspond to the points $q_{0}x-q$ contained in $R_{j}$, by the absolute value of the $q_{0}$ component. If the $q_{0}$ components are equal then order by $q_{1}$ and so on. Suppose that the vector $(m_{0}, \dots ,m_{n})$ is the smallest by our ordering. Then for all other vectors $(r_{0}, \dots , r_{n})$ contained in $R_{j}$ we have that
\begin{align*}
|k_{j,i}-(m_{0}x_{i}-m_{i})-(k_{j,i}-(r_{0}x_{i}-r_{i})|_{p}&<p^{-t_{i}}, \\
|(r_{0}-m_{0})x_{i}-(r_{i}-m_{i})|_{p}&<p^{-t_{i}} \leq N^{-\tau_{i}}.
\end{align*}
Hence the vectors $(r_{0}-m_{0}, \dots , r_{n}-m_{n}) \in \Z^{n+1}$ solve the inequality of $\Qc(\bx, \bt, N)$. Further $(r_{i}-m_{i}) \in [-N,N]$, and by the ordering stated above $r_{0}-m_{0} \in [0,N]$. To exclude the case where $r_{0}-m_{0}=0$ observe that each $\tau_{i}>1$ and so we would have that
\begin{equation*}
N^{-1} \leq |r_{i}-m_{i}|_{p} <p^{-t_{i}} <N^{-1}
\end{equation*}
for $1 \leq i \leq n$, a contradiction. The above argument yields $p^{-n}N^{n+1-\sum_{i=1}^{n}\tau_{i}}-1$ such points, completing the proof.
\qed \\
\vspace{1ex} \par
Lemma \ref{lemma1} is also a relatively simple proof. We are unable to find a proof that uses a similar argument, however we suspect such style of result has been used before. \\
\vspace{1ex} \\
\textit{Proof of Lemma \ref{lemma1}}:
We use a proof by contradiction. Suppose that 
\begin{equation} \label{contra}
\# \Qc(x,\tau,N)>2N^{\tau(x)-\tau +\epsilon}.
\end{equation}
We use the following notations. Let $X \in \N$ be an integer such that 
\begin{equation*}
|x-X|_{p} < p^{-M},
\end{equation*}
for some suitably large $M \in \N$. Define $V_{N}^{+}$ and $V_{N}^{-}$ to be the sets
\begin{align*}
V_{N}^{+}:= \{ (q,q_{1}) \in \N \times \Z: 0<q \leq N, \, 0 \leq q_{1} \leq N, \, \}, \\
V_{N}^{-}:= \{ (q,q_{1}) \in \N \times \Z: 0<q \leq N, \, -N \leq q_{1} \leq 0, \, \}. 
\end{align*}
Let $t \in \N$ be the integer such that
\begin{equation*} 
p^{-t} \leq N^{-\tau} < p^{-t+1},
\end{equation*}
and similarly $k \in \N$ be the integer such that
\begin{equation*} 
p^{-k} \leq N^{-(\tau(x)+\epsilon)} < p^{-k+1}.
\end{equation*}
Note that as $\tau(x)>\tau$, we have that $k \geq t$, and so $p^{k-t} \in \N$. Further, observe that
\begin{equation} \label{k-t}
p^{k-t} < pN^{\tau(x)-\tau+\epsilon}.
\end{equation}
Lastly, by the definition of $\tau(x)$, we have that there exists finitely many $Q \in \N$ such that 
\begin{equation*}
|qx-q_{1}|_{p}< Q^{-(\tau(x)+ \epsilon)},
\end{equation*}
for $0 <q, |q_{1}| \leq Q$. Hence our 'sufficiently large $N_{0}$' is the value of $N_{0}$ such that for any pair $0 <q,|q_{1}| \leq N$,
\begin{equation} \label{bound}
|qx-q_{1}|_{p} \geq N^{-(\tau(x)+\epsilon)},
\end{equation}
for all $\epsilon>0$. Consider the set of points in $\Qc(x, \tau, N)$. Note that $(q,q_{1}) \in \Qc(x, \tau, N)$ if and only if $(q,q_{1}) \in V_{N}^{+} \cup V_{N}^{-}$, and
\begin{equation} \label{congruence}
qX-q_{1} \equiv 0 \mod p^{t}.
\end{equation}
Thus, for all $(q,q_{1}) \in \Qc(x,\tau,N)$ we have that
\begin{equation*}
qX-q_{1}=\lambda p^{t},
\end{equation*}
for some $\lambda \in \Z$. Split the set of points in $\Qc(x,\tau,N)$ into two disjoint sets, the set of pairs in $V_{N}^{+}$, and the set of pairs in $V_{N}^{-}$. As there are greater than $2N^{\tau(x)-\tau+\epsilon}$ pairs, at least one of the sets has greater than $N^{\tau(x)-\tau+\epsilon}$ pairs. Without loss of generality assume such set of points belong in $V_{N}^{+}$. Considering the range of values of $\lambda p^{t}$ there are $p^{k-t}$ possible values of $\lambda p^{t}$ modulo $p^{k}$. By \eqref{contra} and \eqref{k-t} we have, by the Pigeon-hole principle, that there exists at least two pairs, say $(a,a_{1})$ and $(b,b_{1})$, such that 
\begin{equation*}
(a-b)X-(a_{1}-b_{1}) \equiv 0 \mod p^{k}.
\end{equation*}
This is equivalent to 
\begin{equation*}
 |(a-b)x-(a_{1}-b_{1})|_{p} \leq p^{-k} \leq N^{-(\tau(x)+\epsilon)},
 \end{equation*}
 with $(a-b,a_{1}-b_{1}) \in V_{N}^{+} \cup V_{N}^{-}$, as $0<a-b \leq N$ by our choice of ordering of $a,b$, and $|a_{1}-b_{1}| \leq N$ by the fact that the pairs $(a,a_{1}),(b,b_{1}) \in V_{N}^{+}$. 
 However, such result contradicts \eqref{bound} which follows from the definition of $\tau(x)$, thus \eqref{contra} must be false.
\qed \\ 
 \subsection{$p$-adic approximation lattices} \label{p-adic_lattice}
Prior to the proof of Theorem \ref{theorem1} we recall some basic definitions and results of Lattice theory that will be needed. Define a lattice $\Lambda$ as a discrete additive subgroup of $\R^{n}$. If $\Lambda \subseteq \Z^{n}$ the $\Lambda$ is an integer lattice. A set of linearly independent vectors $b_{1}, \dots, b_{n}$ that generate $\Lambda$ is called a basis of $\Lambda$. Let $B$ be a $n\times n$ matrix with columns $b_{i}$, then call $B$ a basis matrix. Define the fundamental region as
\begin{equation*}
\f(B):= \left\{ \sum_{i=1}^{n} a_{i}b_{i} : a_{i} \in \R, \, \, 0 \leq a_{i} < 1 \right\}.
\end{equation*}
A standard result of Lattice theory states that if $B$ is a basis matrix for $\Lambda$ then $\f(B)$ contains no lattice points other than the origin (see Chapter 3, Lemma 6 of \cite{C12}). \par 
 The volume of the fundamental region can be found by taking the determinant of the basis matrix, that is $vol(\f(B))=|\det B|$. A basis matrix is not unique for each $\Lambda$, however for any lattice $\Lambda$ the volume of the fundamental region is the same regardless of choice of basis matrix. For this reason the notation $vol(\f(B))=|\det \Lambda |$ is used. If $U \in \Z^{n\times n}$ is a unimodular matrix and $B_{1}$ is a basis matrix for $\Lambda$ then $B_{2}=B_{1}U$ is also a basis matrix for $\Lambda$.\par
 One property of lattices that are incredibly useful are the successive minima of a lattice. Let $B_{n}=B(0,1)$ denote the $n$-dimensional unit sphere. For $c \in \R_{+}$ we use the notation $cB_{n}=B(0,c)$. Define the successive minima of a lattice $\Lambda \subset \R^{n}$ of rank $n$ as the set of values
\begin{equation*}
\lambda_{i}(\Lambda):=\min \{ \lambda >0: \dim( \Lambda \cap \lambda B) \geq i \},
\end{equation*}
for $i=1, \dots, n$. By Minkowski's inequalities on the successive minima (see e.g. \cite{HF93}) we have that
\begin{equation} \label{mink_bounds}
vol(B_{n})\prod_{i=1}^{n} \lambda_{i}(\Lambda) \leq 2^{n} \det \Lambda.
\end{equation}
For a count on the number of lattice points within a convex body we have the follow theorem due to Blichfeldt \cite{B19}.
\begin{theorem} \label{blichfeldt}
 Let $\Lambda \subset \R^{n}$ be a lattice of rank $n$ and let $V\subset \R^{n}$ be a convex body such that $rank(\Lambda \cap V)=n$.  Then
\begin{equation*}
\# (\Lambda \cap V) \leq n! \frac{vol(V)}{\det \Lambda}+n.
\end{equation*}
\end{theorem}
The constant for such estimate can be excessively large, however in our use of the Theorem the size of such constant is irrelevant. \par
In 1993 an alternative lattice counting theorem was proven by Betke, Henk and Wills \cite{BHW93}, which utilised the properties of the successive minima. This result was further generalised by Henk \cite{H02}, giving us the following theorem.
\begin{theorem} \label{henk_count}
Let $n\geq 2$, $B(0,K)$ a $n$-dimensional ball of radius $K>0$ centred at the origin and $\Lambda$ a $n$-dimensional lattice. Then
\begin{equation*}
\#(\Lambda \cap K) < 2^{n-1}\prod_{i=1}^{n}\left\lfloor \frac{2K}{\lambda_{i}(\Lambda)}+1 \right\rfloor.
\end{equation*}
\end{theorem}
We remark that if $rank(\Lambda \cap B(0,M))<n$ then we must have at least that $\lambda_{n}(\Lambda) \geq M$. Thus the $n^{th}$ value of the product in Theorem \ref{henk_count} would be bounded above by $3$, a point we make use of later on. \par 
For the proof of Theorem \ref{theorem1} we use $p$-adic approximation lattices. First discovered by de Weger \cite{W86} who used them to prove a variety of results in classical $p$-adic Diophantine approximation, including the $p$-adic analogue of Hurwitz Theorem. Recently $n$-dimensional forms of $p$-adic approximation lattices have been used to provided lattice based cryptosystems \cite{IKN17, IN17}. In these papers both dual and simultaneous approximation lattices were discussed. In particular Dirichlet-style exponents were proven for simultaneous and dual approximation. \par
For a $n$-tuple of approximation functions $\Psi=(\psi_{1}, \dots, \psi_{n})$, an integer $N \in \N$, and a fixed $\bx=(x_{1}, \dots, x_{n}) \in \Zp^{n}$ define the $\Psi$-approximation lattice $\Lambda_{N, \bx}$ by
\begin{equation*}
\Lambda_{N, \bx}= \{ (a_{0}, \dots, a_{n}) \in \Z^{n+1}: |a_{0}x_{i}-a_{i}|_{p} \leq \psi_{i}(N), \, \, 1 \leq i \leq n \}.
\end{equation*}
Observe that 
\begin{equation*}
\Qc(\bx,\Psi, N) \subseteq \Lambda_{N,\bx} \cap B(0,\sqrt{n}N),
\end{equation*}
since the euclidean ball $B(0, \sqrt{n}N)$ contains all integer points satisfying $\max_{0 \leq i \leq n}|q_{i}| \leq N$. \par 
 For any $\bx \in \Zp^{n}$ we may write each $x_{j}$ as the $p$-adic expansion
\begin{equation*}
x_{j}=\sum_{i=0}^{\infty}x_{j,i}p^{i}, \quad x_{j,i} \in \{ 0, 1, \dots, p-1 \}.
\end{equation*}
Let $X_{j,N} \in \Z$ be the integer
\begin{equation*}
X_{j,N}=\sum_{i=0}^{t_{j}}x_{j,i}p^{i},
\end{equation*}
where each $t_{j} \in \N$ is the unique value associated with $N$ satisfying
\begin{equation} \label{psi_t}
p^{-t_{j}} \leq \psi_{j}(N) < p^{-t_{j}+1}.
\end{equation}
Lastly, for each $1 \leq j \leq n$ let $\psi_{j,N}^{*}=p^{t_{j}}$. Then the set of vectors
\begin{equation} \label{basis}
B=\left\{ \left( \begin{array}{c}
1 \\
X_{1,N} \\
\vdots \\
X_{n,N}
\end{array} \right),
\left( \begin{array}{c}
0 \\
\psi_{1,N}^{*} \\
\vdots \\
0
\end{array} \right), \dots , 
\left( \begin{array}{c}
0 \\
0 \\
\vdots \\
\psi_{n,N}^{*}
\end{array} \right) \right\},
\end{equation}
form a basis for $\Lambda_{N, \bx}$. The set $B$ can be proven to be a basis by considering the fundamental region $\F(B)$ and showing the only lattice point contained is $\textbf{0}$. Given such basis we can deduce that
\begin{equation*}
|\det \Lambda_{N, \bx} |=\prod_{i=1}^{n}\psi_{i,N}^{*} \asymp \left(\prod_{i=1}^{n} \psi_{i}(N) \right)^{-1},
\end{equation*}
where the implied constants can be easily found using \eqref{psi_t} to obtain
\begin{equation} \label{det_Lambda}
\left(\prod_{i=1}^{n} \psi_{i}(M)\right)^{-1}\leq |\det \Lambda_{N, \bx}| \leq p^{n}\left(\prod_{i=1}^{n} \psi_{i}(N)\right)^{-1}.
\end{equation}
In the simultaneous case, $\Psi=(\psi, \dots, \psi)$, it was proven in \cite{IN17} that 
\begin{equation*}
\lambda_{1}(\Lambda_{N, \bx}) \ll \psi(N)^{-\frac{n}{n+1}}.
\end{equation*}
In order to prove Theorem \ref{theorem1} we find a lower bound on $\lambda_{1}(\Lambda_{N, \bx})$ by considering $\bx \in \Zp$ satisfying certain Diophantine exponent properties. For completeness we also prove the upper bound.

\begin{lemma} \label{bounds_Lmin}
Let $\Lambda_{N, \bx}$ be defined above with $\tau(\bx)=n+1$, and suppose that
\begin{equation*}
\prod_{i=1}^{n}\psi_{i}(N) < N^{-n}.
\end{equation*}
 Then for any $\varepsilon>0$ the exists sufficiently large $N_{0} \in \N$ such that for all $N \geq N_{0}$,
\begin{equation*}
\left( \frac{1}{\prod_{i=1}^{n}\psi_{i}(N)} \right)^{\frac{1}{n+1}-\varepsilon} \leq \lambda_{1}(\Lambda_{N, \bx}) \leq C_{2}\left( \frac{1}{\prod_{i=1}^{n}\psi_{i}(N)} \right)^{\frac{1}{n+1}},
\end{equation*}
where 
\begin{equation*}
C_{2}=2\left(\frac{ \Gamma\left(\frac{n+1}{2}+1 \right)p^{n}}{\pi^{\frac{n+1}{2}}}\right)^{\frac{1}{n+1}}.
\end{equation*}
\end{lemma}
As will become clear in the proof below the condition that $\tau(\bx)=n+1$ is only necessary in the lower bound result.
\begin{proof}
We prove the upper bound case first. Such proof is a standard application of Minkowski's first Theorem on successive minima and follows almost immediately by the above calculation of $\det(\Lambda_{M, \bx})$. Concisely, we have that
\begin{equation*}
\lambda_{1}(\Lambda_{N,\bx})^{n+1}vol(B(0,1)) \leq 2^{n+1} \det(\Lambda_{N, \bx}).
\end{equation*}
 Rearranging for $\lambda_{1}(\Lambda_{N,\bx})$, using \eqref{det_Lambda}, and recalling the volume of an $n+1$-ball we obtain our result. \par
 We now prove the lower bound. For any $\bx \in \Zp^{n}$ observe that
\begin{equation*}
\prod_{i=1}^{n}|q_{0}x_{i}-q_{i}|_{p} < M^{-(n+1)},
\end{equation*}
for infinitely many $M$ (see for example Lemma \ref{mink5}). Further, since $\tau(\bx)=n+1$ there exists $N_{0}$ such that for all $N \geq N_{0}$ then any rational integer vectors $(q_{0},\dots , q_{n})$ satisfying $\max_{0 \leq i \leq n}|q_{i}| \leq N$ we have that
\begin{equation} \label{extremal_contradiction}
\prod_{i=1}^{n}|q_{0}x_{i}-q_{i}|_{p} \geq N^{-(n+1+\varepsilon)}
\end{equation}
for some $\varepsilon>0$. Choose $N$ sufficiently large such that
\begin{equation*}
N_{0} \leq \left(\prod_{i=1}^{n} \psi_{i}(N)\right)^{-\left(\frac{1}{n+1}-\varepsilon\right)}.
\end{equation*}
 Such $N$ is possible since $\prod_{i=1}^{n} \psi_{i}(N)<N^{-n}$ and so the value on the RHS of the above inequality tends to infinity as $N \to \infty$ for any small $\varepsilon$ ($\varepsilon<\frac{1}{n(n+1)}$). \par
Suppose that $(q_{0}, \dots , q_{n})$ is a minimum length vector of $\Lambda_{N,\bx}$, then note that $\lambda_{1}(\Lambda_{N, \bx}) \geq \max_{1 \leq i \leq n}|q_{i}|$ due to the euclidean nature of $\lambda_{1}(\Lambda_{N,\bx})$. Suppose that
\begin{equation} \label{false}
\max_{1 \leq i \leq n}|q_{i}| < \left(\prod_{i=1}^{n} \psi_{i}(N)\right)^{-\left(\frac{1}{n+1}-\varepsilon\right)}.
\end{equation}
We prove \eqref{false} to be false. Observe that
\begin{equation*}
\prod_{i=1}^{n}|q_{0}x_{i}-q_{i}|_{p}<\prod_{i=1}^{n}\psi_{i}(N),
\end{equation*}
since $(q_{0},\dots , q_{n}) \in \Lambda_{N, \bx}$. Then
\begin{equation*}
\prod_{i=1}^{n}|q_{0}x_{i}-q_{i}|_{p}<\left(\left(\prod_{i=1}^{n}\psi_{i}(N)\right)^{-\left( \frac{1}{n+1}- \varepsilon \right)}\right)^{-\frac{n+1}{1-\varepsilon(n+1)}}.
\end{equation*}
But this contradicts \eqref{extremal_contradiction}. So we must have that \eqref{false} is false, and so
\begin{equation*}
\lambda_{1}(\Lambda_{N,\bx}) \geq \left(\prod_{i=1}^{n} \psi_{i}(N)\right)^{-\left(\frac{1}{n+1}-\varepsilon\right)},
\end{equation*}
completing the proof. 
\end{proof}
Given Lemma \ref{bounds_Lmin} we can proceed with the following. \\
 \vspace{1ex}
\textit{Proof of Theorem \ref{theorem1}}:
   For $N \geq N_{0}$, where $N_{0}$ is chosen by Lemma \ref{bounds_Lmin}, consider the following two cases:
\begin{enumerate}[i)]
\item $rank(\Lambda_{N,\bx} \cap B(0,\sqrt{n}N))=n+1$: By Theorem~\ref{blichfeldt}  we have that
\begin{align*}
\# (\Lambda_{N,\bx} \cap B(0,\sqrt{n}N)) & \leq (n+1)! \frac{vol(B(0,\sqrt{n}N))}{\det \Lambda_{N,\bx}}+n+1, \\
& \leq \frac{(n+1)! \pi^{n/2}\sqrt{n}^{n+1}}{\Gamma\left(\frac{n}{2}+1 \right)}N^{n+1}.(\prod_{i=1}^{n}\psi_{i,N}^{*})^{-1}+n+1, \\
& \leq  \frac{(n+2)! \pi^{n/2}\sqrt{n}^{n+1}}{\Gamma\left(\frac{n}{2}+1 \right)} N^{n+1}\prod_{i=1}^{n}\psi_{i}(N).
\end{align*}
Note that the last inequality follows since $\prod_{i=1}^{n}\psi_{i}(N) > N^{-(n+1-\varepsilon)}$. This proves Theorem~\ref{theorem1} for the rank $n+1$ case. \\
\item $rank(\Lambda_{N,\bx} \cap B(0,\sqrt{n}N))<n+1$: Since $rank(\Lambda_{N,\bx} \cap B(0,\sqrt{n}N))<n+1$ we must have $\lambda_{n+1}(\Lambda_{N,\bx})>\sqrt{n}N$. Hence, by the remark made previously, the final product on the right of Theorem \ref{henk_count} is less than or equal to $3$. Furthermore, for each $\lambda_{i}(\Lambda_{N,\bx})$, $1 \leq i \leq n$ we have that
\begin{align*}
\lambda_{n}(\Lambda_{N,\bx}) \geq \dots \geq \lambda_{1}(\Lambda_{N,\bx}) & \overset{(a)}{\geq} \left( \frac{1}{\prod_{i=1}^{n}\psi_{i}(N)} \right)^{\frac{1}{n+1}-\varepsilon}, \\
& \overset{(b)}{\geq} \left( \frac{1}{N \prod_{i=1}^{n}\psi_{i}(N)} \right)^{1/n},
\end{align*}
where $(a)$ follows from Lemma \ref{bounds_Lmin} and $(b)$ follows since
\begin{align*}
 \left( \frac{1}{\prod_{i=1}^{n}\psi_{i}(N)} \right)^{\frac{1}{n+1}-\varepsilon} & \geq \left( \left(\frac{1}{\prod_{i=1}^{n}\psi_{i}(N)}\right)^{n-\varepsilon n(n+1)} \right)^{\frac{1}{n(n+1)}}, \\ & \geq \left( \frac{1}{N^{n+1}\prod_{i=1}^{n}\psi_{i}(N)}\left( \frac{1}{\prod_{i=1}^{n}\psi_{i}(N)} \right)^{n}\right)^{\frac{1}{n(n+1)}}.
 \end{align*}
 combining the two ideas above, and Theorem \ref{henk_count}, we have that
 \begin{align*}
 \# (\Lambda_{N, \bx} \cap B(0,\sqrt{n}N)) & < 2^{n}3 \prod_{i=1}^{n}\left( \frac{2\sqrt{n}N}{\lambda_{1}(\Lambda_{N,\bx})}+1 \right), \\
 &< 2^{n}3 \left(2\sqrt{n}N^{1+1/n}\left(\prod_{i=1}^{n}\psi_{i}(N) \right)^{1/n}+1 \right)^{n}, \\
 &<3(6\sqrt{n})^{n}N^{n+1}\prod_{i=1}^{n}\psi_{i}(N).
 \end{align*}
 
  \end{enumerate}
Thus, in either case $i)$ or $ii)$ we have that
 \begin{equation*}
  \# (\Lambda_{N, \bx} \cap B(0,\sqrt{n}N)) \leq C_{1} N^{n+1}\prod_{i=1}^{n}\psi_{i}(N),
  \end{equation*}
  with
  \begin{equation*}
  C_{1}=\max\left\{ 3(6\sqrt{n})^{n}, \frac{(n+2)! \pi^{n/2}\sqrt{n}^{n+1}}{\Gamma\left(\frac{n}{2}+1 \right)}\right\}.
  \end{equation*}
  
 \qed

\section{Concluding remarks}
This article provides sharp bounds on the number of rational points close to almost all $p$-adic integers. While this result allows us to find simultaneous $p$-adic Diophantine approximation results on coordinate hyperplanes, it falls a long way short of providing results for Diophantine approximation sets on curves and manifolds. It is hoped the techniques used in this paper could be used to find rational points close to manifolds, we intend to follow this up with a subsequent paper. \par 

\bibliographystyle{plain}
\bibliography{biblio}

\end{document}

%% file: header.tex
\usepackage{graphicx}
\usepackage{amsmath}
\usepackage{amsfonts}
\usepackage{amsthm}
\usepackage{hyperref} 
\usepackage{amsopn,amssymb,url,mathrsfs,a4wide}
\usepackage{tikz}
\usepackage{booktabs}
\usepackage[margin=0.8in]{geometry}
\usepackage{parskip}
\usepackage{fancyvrb}
\usepackage{listings}

\newcounter{ctr}

\newcounter{ctr1}

\newcounter{ctr2}

\newcounter{ctr3}

\newtheorem{definition}{Definition}[section] 
\newtheorem{theorem}[definition]{Theorem}
\newenvironment{theorem*}[1]{{\bf Theorem #1} \begin{itshape}}{\end{itshape}}
\newtheorem{lemma}[definition]{Lemma}
\newtheorem{corollary}[definition]{Corollary}
\newenvironment{corollary*}[1]{{\bf Corollary #1} \begin{itshape}}{\end{itshape}}
\newtheorem{proposition}[definition]{Proposition}
\newenvironment{proposition*}[1]{{\bf Proposition #1} \begin{itshape}}{\end{itshape}}
\newtheorem{remark}[definition]{Remark}

\newcommand{\ud}{\, {\rm d} \kern-.015em }

\newcommand{\modulus}[1]{\left| \kern.05em #1 \kern.05em \right|}
\newcommand{\norm}[1]{\left\| \kern.05em #1 \kern.05em \right\|}
\newcommand{\inner}[1]{\left\langle \kern.05em #1 \kern.05em \right\rangle }

\newcommand{\pick}[2]{\renewcommand{\arraystretch}{0.6}
\left( \kern-.4em \begin{array}{c} #1 \\ #2 \end{array} \kern-.4em \right) }

\newsavebox{\FVerbBox}



%% file: help.tex
\newcommand{\ha}{\mathcal{H}}

\newcommand{\Q}{\mathbb{Q}}

\newcommand{\f}{\mathcal{F}}
\newcommand{\R}{\mathbb{R}}
\newcommand{\N}{\mathbb{N}}
\newcommand{\Z}{\mathbb{Z}}
\newcommand{\bx}{\boldsymbol{x}}
\newcommand{\bt}{\boldsymbol{\tau}}

\newcommand{\ba}{\boldsymbol{a}}
\newcommand{\Cf}{\mathcal{C}_{f}}

\newcommand{\Zp}{\mathbb{Z}_{p}}
\newcommand{\Qp}{\mathbb{Q}_{p}}
\newcommand{\W}{\mathcal{W}}
\newcommand{\A}{\mathcal{A}}
\newcommand{\M}{\mathcal{M}}

\newcommand{\Ap}{\mathfrak{A}}

\newcommand{\Wp}{\mathfrak{W}}

\newcommand{\F}{\mathcal{F}}

\newcommand{\bq}{\boldsymbol{q}}

\newcommand{\bal}{\boldsymbol{\alpha}}
\newcommand{\tv}{\tilde{v}}